\documentclass[a4paper,12pt]{article}
\usepackage[english]{babel}
\usepackage[T1]{fontenc}
\usepackage[latin1]{inputenc}
\usepackage{graphics, enumerate}
\usepackage{amsmath,amsfonts,amstext,amssymb,color,epsfig,bbm}

\numberwithin{equation}{section}

\newtheorem{stheorem}{Theorem}[section]
\newtheorem{definition}[stheorem]{Definition}
\newtheorem{proposition}[stheorem]{Proposition}
\newtheorem{lemma}[stheorem]{Lemma}
\newtheorem{corollary}[stheorem]{Corollary}
\newenvironment{proof}{{\it \underline{proof}}\ :\ }{\qed $\\$}

\newenvironment{remark}{{\ \ \ \it Remark : }}{}
\def\P {\mathbb {P}}

\def\pp#1{\overset{..}#1}
\def\phi{\varphi}
\def\D{\Delta}

\def\tpar{/ \! \! /}

\def\R{\mathbb{R}}
\def\E{\mathbb{E}}
\def\N{\mathbb{N}}
\def\1{\mathbbm{1}}

\def\qed{\hbox{$\vcenter{\vbox{
   \hrule height 0.4pt\hbox{\vrule width 0.4pt height 6pt
    \kern5pt\vrule width 0.4pt}\hrule height 0.4pt}}$}}

\title{Some stochastic process without birth, linked to the mean curvature flow}
\author{A.K. Coulibaly}

\begin{document}
\date{}
\maketitle



\begin{abstract}
Using Huisken results about the mean curvature flow on a strictly convex hypersurface, and Kendall-Cranston coupling, we will build a stochastic process without birth, and show that there exists a unique law of such process. This process has many similarities with the circular Brownian motions studied by \'Emery, Schachermayer, and Arnaudon. In general, this process is not a stationary process, it is linked with some differential equation without initial condition. We will show that this differential equation has a unique solution up to a multiplicative constant.
\end{abstract}

\section{Tools and first properties}
Let $M$ be a Riemannian compact $n$-manifold without boundary, which is smoothly embedded in $\R^{n+1}$, and $n \ge 2$. Denote by $F_{0} $ the embedding function:
\begin{equation*}
 F_{0}: M \hookrightarrow \R^{n+1}.
\end{equation*}
Consider the flow defined by:
\begin{equation} \label{mean-eq}
\left\{ \begin{array}{l} 
 \partial_{t}F(t,x) = -H_{\nu}(t,x) \overset{\rightarrow}{\nu}(t,x) \\
 F(0,x) = F_{0}(x) .
         \end{array} 
\right.
 \end{equation}
 Let $ M_{t} = F(t,M)$, we identify $M$ with $M_{0}$ and $F_{0}$ with $Id$. In the previous equation (\ref{mean-eq}), ${\nu}(t,x) $ is the outer unit normal at $F(t,x)$ on $ M_{t} $, and  $H_{\nu}(t,x) $ is the mean curvature at $F(t,x)$ on $ M_{t} $ in the direction ${\nu}(t,x) $, (i.e. $H_{\nu}(x) = \text{ trace }(S_{\nu}(x))$ where $S_{\nu}$ is the second fundamental form, for definition see \cite{Lee}).

\begin{remark}
 In this paper we take this point of view of the mean curvature flow (see \cite{Huisken} for existence, and related result).
 Many other authors give a different point of view for this equation. The viscosity solution (see \cite{evsp1},\cite{evsp2},\cite{evsp3},\cite{evsp4},\cite{evsonsou}) generalizes the solution after the explosion time and gives a uniqueness solution which is also contained in Brakke family of solutions and passes the singularity. We will just look at the smooth solution until the explosion time.
\end{remark}
\\

%
 As usual we call $ M_{t}$ the motion by mean curvature. For self-completeness, we include a proof of the next lemma, although it is well-known.

\begin{lemma}
 Let $(M,g)$ be a Riemannian manifold isometrically embedded in $\R^{n+1}$. Denote $\iota $ the isometry: 
       $$ (M , g) \overset {\iota} {\hookrightarrow} \R^{n+1}.$$ 
then:
 \begin{equation}  \label{lap}  
\forall x \in M ,\quad  \Delta \iota (x) = - H_{\nu}(x) \vec {\nu} (x).
 \end{equation}
Where $\Delta $ is the Laplace-Beltrami operator associated to the metric $g$.
\end{lemma}

\begin{proof}
 By the flatness of target manifold, we have
$$ \D \iota (x) = \left(\begin{array}{c}
                        \D \iota^1 (x) \\
                         \vdots \\
                         \D \iota^{n+1} (x) \\
                        \end{array} \right)  $$ 
and 
$\begin{array}{rcl}
\D \iota^j (x)  &=& \sum_{i=1}^{n} \frac{d}{dt^2}\big|_{t=0} \iota^j (\gamma_i(t)),\\
\end{array}$
 where $\gamma_i(t)$ is a geodesic in $M$ such that $\gamma_i(0)= x$ and $\dot\gamma_i(0)=A_{i} $ and $A_{i}$ is a orthogonal basis of $ T_{x}M$.
 By definition of a geodesic we obtain:
$$\D \iota (x) \perp T_{\iota(x)}(\iota(M)), $$
so there exists a function $\beta$ such that $\D \iota (x) = \beta(x) \vec{\nu}(x) $. We compute $\beta$ as follows:
$$\begin{array}{rcl}
\beta (x) &=& \langle   \D \iota(x) , \vec{\nu}(x)\rangle \\
          &=& \sum_{i=1}^{n} \langle \frac{d}{dt^2}\big|_{t=0} \iota (\gamma_i(t)) , \vec{\nu}(x)\rangle \\
          &=& \sum_{i=1}^{n} \langle \nabla^{\R^n}_{\dot{\iota(\gamma_i(t))}}\dot{\iota(\gamma_i(t))} \big|_{t=0}, \vec{\nu}(x)\rangle \\
          &=& \sum_{i=1}^{n} - \langle \dot{\iota(\gamma_i(t))} , \nabla^{\R^n}_{\dot{\iota(\gamma_i(t))}}  \vec{\nu} \rangle \big|_{t=0} \text{,   metric connection}\\
          &=& \sum_{i=1}^{n} - \langle \dot{\iota(\gamma_i(t))} , (\nabla^{\R^n}_{\dot{\iota(\gamma_i(t))}}  \vec{\nu})^{\top} \rangle \big|_{t=0} \\
          &=& - \text{ trace }(S_{\nu}(x)). \\
\end{array}$$
\end{proof}

To give a parabolic interpretation of this equation (\ref{mean-eq}), let us define a family of metrics $ g(t) $ on $M$ which is the pull-back by $ F(t,.)$ of the induced metric on $ M_{t}$. Using the previous lemma we rewrite the equation as in (\cite{Huisken}):

\begin{equation}
\left\{ \begin{array}{l} 
 \partial_{t}F(t,x) = \D_{t} F(t,x)\\
 F(0,x) = F_{0}(x) 
         \end{array} 
\right.
 \end{equation} 
where $\D_{t} $ is the Laplace-Beltrami operator associated to the metric $g(t)$.

\begin{remark}
Sometimes we will use a probabilistic convention, consisting in putting $\frac12$ before the Laplacian (which just changes the time and makes the calculus more synthetic), sometimes we will use geometric convention.
\end{remark}


We call $T_{c}$ the explosion time of the mean curvature flow, let $T < T_{c}$, and $g(t)$ be the family of metrics defined above.
Let $ (W^{i})_{1 \le i \le n }$ be a $\R^n$-valued Brownian motion. We recall from \cite{metric} the definition of the $g(t)$-Brownian motion in $M$ started at $x$, denoted by $g(t)$-BM($x$):

\begin{definition}
Let us take a filtered probability space $(\Omega, (\mathcal{F}_{t})_{t \ge 0},\mathcal{F},\mathbb{P})$ and a $C^{1,2}$-family $g(t)_{t\in [0,T[}$ of metrics over $M$. A $M$-valued process $X(x)$ defined on $ \Omega \times [0,T[$ is called a $g(t)$ Brownian motion in $M$ started at $x\in M$ if $X(x)$ is continuous, adapted and for every smooth function $f$,
  $$ f(X_{s}(x)) -f(x) - \frac12 \int_{0}^{s} \D_{t}f (X_{t}(x)) dt  $$  
is a local martingale vanishing at $0$.  
\end{definition}    
 
We give a proposition which yields a characterization of mean curvature flow by the $g(t)$ Brownian motion.

\begin{proposition}\label{prop-pullbac} 
 Let $M$ be an n-dimensional manifold isometrically embedded in $\R ^{n+1}$. Consider the application:
 $$ F: [0,T[ \times M  \rightarrow  \R^{n+1} $$
such that $F(t,.)$ are diffeomorphisms, and the family of metrics $g(t)$ over $M$, which is the pull-back by $ F(t,.)$ of the induced metric on $ M_{t}=F(t,M)$. Then the following items are equivalent:
\begin{enumerate}[i)]
\item  $F(t , .)$ is a solution of mean curvature flow
\item $\forall x_0 \in M $, $\forall T \in [0 , T_{c}[ $, let $\tilde{g}_t^{T} = \frac12 g_{T - t}$  and $X^{T}(x_0) $ be a   
$(\tilde{g}_t^{T})_{t \in [0, T]}\text{-BM}(x_0)$, then: 
  $$ Y^{T}_t = F(T -t , X^{T}_t(x_0)) $$ 
is a local martingale in $\R ^{n+1} $. 
\end{enumerate}
\end{proposition}
 
\begin{proof}
 By definition we have a sequence of isometries:
$$ F(t , .): (M , g_t) \tilde{\rightarrow} M_t \hookrightarrow \R^{n +1} $$ 
Let $x_{0} \in M$ and $T \in [0,T_{c}[$ and  $X^{T}(x_0)$ a $ (\tilde{g}_t^{T})_{t \in [0, T]}\text{-BM}(x_0). $ 
We just compute the It\^o differential of:
 $$ Y^{T , i}_t = F^{i}(T -t , X^{T }_t(x_0)), $$
that is to say:
$$ \begin{array} {lcl}
d(Y^{T , i}_t ) &=&  - \frac{\partial }{\partial t} F^{i}(T - t , X^{T }_t(x_0) ) dt + d (F_{T -t}^{i} ( X^{T }_t(x_0)) \\
                      & \underset {d \cal{M}}{\equiv }&  - \frac{\partial }{\partial t} F^{i}(T - t , X^{T }_t(x_0) ) dt  + \frac12 \D_{\tilde{g}_t} F_{T -t}^{i} (X^{T }_t(x_0) ) dt  \\
                      & \underset {d \cal{M}}{\equiv }&  - \frac{\partial }{\partial t} F^{i}(T - t , X^{T }_t(x_0) ) dt  +  \D_{g_{T - t} }F_{T -t}^{i} (X^{T }_t(x_0) ) dt\\
                      & \underset {d \cal{M}}{\equiv }&  0.\\
\end{array}\\$$
Therefore $ Y^{T}_t $ is a local martingale.

Let us show the converse. Let $x_{0} \in M$ and $T \in [0,T_{c}[$ and  $X^{T}(x_0)$ a $ (\tilde{g}_t^{T})_{t \in [0, T]}\text{-BM}(x_0) $, $Y^{T , i}_t $ is a local martingale so almost surely, for all $t \in [0,T]$:
$$  - \frac{\partial }{\partial t} F^{i}(T - t , X^{T }_t(x_0) ) dt  +  \D_{g_{T - t} }F_{T -t}^{i} (X^{T }_t(x_0) ) dt=0$$
so that for all $s \in [0,T]$, by integrating we get  
$$ \int_0^s - \frac{\partial }{\partial t} F^{i}(T - t , X^{T }_t(x_0) ) dt  +  \D_{g_{T - t} }F_{T -t}^{i} (X^{T }_t(x_0) ) dt =0 $$
the continuity of every $g(t)$-Brownian motion yields,
 $$ - \frac{\partial }{\partial t} F^{i}(T  , x_0 )   +  \D_{g_{T } }F_{T }^{i}( x_0 ) = 0 .$$

 \end{proof}
In order to apply this proposition, we give an estimation of the explosion time. It is also a consequence of a maximum principle, which is explicitly contained in the $g(t)$-Brownian Motion.

The quadratic covariation of $Y^{T}_t $ is given by:

\begin{proposition}\label{crochet}
Let  $Y^{T}_t $ be defined as before, then the quadratic covariation of  $Y^{T}_t $ for the usual scalar product in $\R^{n+1}$ is:
$$\langle d Y^{T}_t ,d Y^{T}_t \rangle = 2n \1_{[0,T]}(t) dt $$ 
\end{proposition}
\begin{proof}
 Let $\tpar^{T}_{0,t} $ be the parallel transport above $X^{T}_{t}$, it is shown in \cite{metric} that it is an isometry :
$$\tpar^{T}_{0,t}: (T_{X_{0}}M, \tilde{g}(0)) \longmapsto  (T_{X_{t}}M, \tilde{g}(t)). $$
Let $(e_{i})_{1 \le  i \le n}$ be a orthonormal basis of $(T_{X_{0}}M,\tilde{g}(0))$, and  $ (W^{i})_{1 \le i \le n }$ be the $\R^n$-valued Brownian motion such that (e.g. \cite{metric}, \cite{Crass-KAT}):
$$*dW_{t} =  \tpar^{T,-1}_{0,t} *d X^{T}_{t},$$ 
and in the It\^o's sense:  
$$dX^{T}_{t} = \tpar^{T}_{0,t}e_{i} dW^{i}_{t}.$$ 
Hence 
$$\begin{array}{lcl}

\langle d Y^{T}_t ,d Y^{T}_t \rangle &=& \langle d (F_{T -t} ( X^{T }_t(x_0))) ,  d (F_{T -t} ( X^{T }_t(x_0))) \rangle \\
&=& \langle d ( X^{T }_t(x_0)) ,  d( X^{T }_t(x_0)) \rangle_{g_{T -t}} \\
&=& \langle d ( X^{T }_t(x_0)) ,  d( X^{T }_t(x_0)) \rangle_{2\tilde{g}_{t }} \\
&=& \langle \sum_{i=1}^{n} \tpar^{T}_{0,t}e_{i} dW^i  , \sum_{j=1}^{n}\tpar^{T}_{0,t} e_{j} dW^j   \rangle_{2\tilde{g}_{t }} \\
&=&\sum_{i=1}^{n}  \langle  \tpar^{T}_{0,t}e_{i} ,\tpar^{T}_{0,t} e_{i}   \rangle_{2\tilde{g}_{t }} dt \\
&=&\sum_{i=1}^{n}  2  dt  \\
&=& 2n dt. \\
\end{array} $$
To go from the first to the second line, we have used the fact that $F_{T -t}$ is a isometry, for the last step we used the isometry of the parallel transport.
\end{proof}

\begin{remark}
 Up to convention we recover the same martingale as in \cite{sontou}.
\end{remark}
\\
An immediate corollary of Proposition \ref{crochet} is the following result, which appears in \cite{Huisken} and \cite{evsp1}. 
\begin{corollary}
 Let $M$ be a compact Riemannian $n$-manifold and $T_{c}$ the explosion time of the mean curvature flow, then:
 $T_{c} \le \frac{diam(M_0)^2}{2n}$ 
\end{corollary}
\begin{proof}
 Recall that the mean curvature flow stays in a compact region, like the smallest ball which contain $M_{0}$, this result is clear in the strictly convex starting manifold and can be found in a general setting using P.L Lions viscosity solution (e.g. theorem 7.1 in \cite{evsp1}).
\\
For all $T \in [0, T_{c}[$ take the previous notation. So by the above recall that:
 $$ \parallel Y^{T}_t  \parallel \leq diam(M_0), $$
then $ Y^{T}_t $ is a true martingale.
And

  $$ \parallel Y^{T}_t  \parallel^2 -   \langle  Y^{T} , Y^{T} \rangle_t  $$ is also a true martingale. Hence: 
$$\E [\parallel Y^{T}_0 \parallel^2] + 2 nT \le diam(M_0)^2 ,$$
we obtain 
$$T \le \frac{diam(M_0)^2}{2n} .$$
\end{proof}

\section{Tightness, and first example on the sphere}
 We now define $(\tilde{g}^{T_{c}})_{t \in ]0,T_{c}]}$-BM in a general setting. When the initial manifold $M_{0}$ is a sphere we use the conformality of the metric, to show that after a deterministic change of time such process is a $] -\infty , T_{c}]$ Brownian motion on the sphere (for existence and definition see \cite{Emery_Schachermayer} and  \cite{Marcfiltration} ). In the next section, we will give a general result of uniqueness when the initial manifold $M_{0}$ is strictly convex.

\begin{definition}
 Let $M$ be an n-dimensional strictly convex manifold (i.e. with a strictly positive definite second fundamental form), $F(t,.)$ the smooth solution of the mean curvature flow, $(M,g(t)) $ the family of metrics constructed by pull-back (as in \ref{prop-pullbac}) and $T_{c}$ the explosion time.
We define a family of processes as follow:\\
$\forall \epsilon \in ]0,T_{c}] $

$$ X_t^{\epsilon}(x_0) = \left\{\begin{array}{l}
                           x_0 \text{  if }  0 < t \le \epsilon \\
\\
                           BM(\epsilon, x_0)_{t} \text{  if }  \epsilon \le t \le T_{c}\\ 

                           \end{array} \right.$$  

where $ BM(\epsilon, x_0)_{t}$ denotes a $\frac12 g(T_{c} - t)$ Brownian motion that starts at $x_{0}$ at time $\epsilon$, and 
\\
$$ Y_t^{\epsilon}(x_0) = \left\{\begin{array}{l}
                           F(T_{c}-\epsilon , x_0) \text{  if } 0 \le t \le \epsilon \\
\\
                           F(T_{c}-t,X_t^{\epsilon}(x_0) ) \text{  if } \epsilon \le t \le T_{c}.\\ 

                           \end{array} \right. \\$$  
\end{definition}
\begin{remark}
We proceed as before because, at the time $T_{c}$, there is no more metric. 
Huisken shows in \cite{Huisken} that in this case:
$$\exists \mathcal{D} \in \R^{n+1} ,\quad s.t.  \quad \forall x_0 \in M ,\qquad \lim_{s \rightarrow T_{c} } F(s , x_0)=\mathcal{D}  $$ 
\end{remark}

\begin{proposition}\label{tight_Y}
 With the same notation as the above definition, there exists at least one martingale $Y^{1}$ in the adherence (for the weak convergence) of $ (Y_.^{\epsilon}(x_0))_{\epsilon>0}$ when $\epsilon$ goes to $0$. Also, every adherence value is a martingale. 
 \end{proposition}
\begin{proof}
  We have:
$$
\left\{ 
\begin{array}{l}
dY_t^{\epsilon}(x_0) = 0 \text{  if t } \le \epsilon \\
                                                     \\
dY_t^{\epsilon}(x_0) = d \mathcal{M} \text{  if t }\ge \epsilon .
\end{array}
\right.
$$
 Where $d \mathcal{M}$ is an It\^o differential of some martingale. This defines a family of martingales. With the same computation as in proposition \ref{crochet}, we get:
$$ \langle d Y^{\epsilon}_t ,d Y^{\epsilon}_t \rangle_{\R^{n+1}} = 2n1_{]\epsilon , T_{c}]}(t) dt \le 2n dt.$$
%
Also by the above remark $Y^{\epsilon}_0$ is tight, hence $ (Y_{.}^{\epsilon}(x_0))_{\epsilon >0}$ is tight. As usual, Prokhorov's theorem implies that one adherence value exists.
We also use Huisken \cite{Huisken} (for the strictly convex manifold) to yield:
\begin{equation}\label{bornee}
  \parallel Y^\epsilon \parallel \le diam(M_0).
\end{equation}
By proposition 1-1 in \cite{JS}  page 481, and the fact that $ (Y^\epsilon)$ are martingales we conclude that all adherence values of 
$ (Y^\epsilon)$ are martingales with respect to the filtration that they generate.
\end{proof}
\\
\begin{remark}
 The above proposition is also valid for arbitrary $M$ that are isometrically embedded in $\R^{n+1}$. Just because the bound \ref{bornee} is also a consequence of theorem 7.1 in \cite{evsp1}.
\end{remark}
\\
We will now derive the tightness of $X^\epsilon_t$ from those of $ (Y^\epsilon)$. 
This purpose will be completed by the next lemma \ref{lem_sko}.
\\
Recall some results of \cite{Huisken}, if $M_{0}$ is a strictly convex manifold then $M_{t}$ is also strictly convex, and  $\forall 0 \le t_{1} < t_{2} < T_{ c}$, $ M_{t_{2}} \subset int(M_{t_{1}})$, where $int$ is the interior of the bounded connected component.
Hence there is a foliation on $\overline{int(M_{0})} $:
$$\bigsqcup_{t \in [0,T_{c}[ } M_{t},$$ 
where $\bigsqcup$ stand for the disjoint union.
\begin{definition}
We note:
$$\mathcal{C}^f(]0,T_{c}],\R^{n+1}) = \{ \gamma \in \mathcal{C}(]0,T_{c}],\R^{n+1}) ,s.t. \quad \gamma(t)\in M_{T_{c}-t} \}. $$

\end{definition}
Noted that $\mathcal{C}^f(]0,T_{c}],\R^n)$ is a closed set of $\mathcal{C}(]0,T_{c}],\R^n)$ for the Skorokhod topology.

\begin{lemma}\label{lem_sko}
 Let $M$ an n-dimensional strictly convex manifold, $F(t,.)$ the smooth solution of the mean curvature flow and $T_{c}$ the explosion time. Then
$$\begin{array}{rlc}
 F: [0,T_{c}[\times M &\longrightarrow&  \bigsqcup_{t \in [0,T_{c}[ }M_{t} \\   
\end{array},$$
is a diffeomorphism in the sense of manifold with boundary.
And,

\begin{center}
$\begin{array}{rlc}
  \Psi: \mathcal{C}^f(]0,T_{c}],\R^n) & \longrightarrow& \mathcal{C}(]0,T_{c}],M ) \\
             \gamma &\longmapsto &  t \mapsto F^{-1}(T_{c}-t , \gamma(t))
\end{array}$
\end{center}
is continuous for the different Skorokhod topologies. To define the Skorokhod topology in $\mathcal{C}(]0,T_{c}],M )$ we could use the initial metric $g(0)$ on $M$.
\end{lemma}

\begin{proof}
 It is clear that $F$ is smooth as a solution of a parabolic equation \cite{Huisken}, and this result has been used above. Its differential is given at each point by:
$$\forall (t,x) \in [0,T_{c}[\times M  , \quad \forall v \in T_{x} M $$
$$DF(t,x) (\frac{\partial}{\partial_{t}} , v) = \frac{\partial}{\partial_{t}} F(t,x) \oplus DF_{t}(x)(v)  $$
where $\frac{\partial}{\partial_{t}} F(t,x) =- H(t,x)\overrightarrow{\nu}(t,x)$, here $\oplus$ stands for + and means that we cannot cancel the sum without cancelling each term. Since there is no ambiguity we write $H(t,x)$ for $H_{\nu}(t,x)$. Recall that $ H(t,x) >0 $.

 For the second part of this lemma, we remark that for $ 0 \le \delta < T_{c}$  
$$F^{-1}: \bigsqcup_{t \in [0,\delta] } M_{t} \longrightarrow [0,\delta]\times M  $$ 
 is Lipschitz (use the bound of the differential on a compact).

Recall also that a sequence converges to a continuous function for Skorokhod topology if and only if it converges to this function locally uniformly.
 We will now show the continuity of $\Psi $. Take a sequence  $ \alpha_{m} $ in $\mathcal{C}^f(]0,T_{c}],\R^{n+1})$ and $\alpha \in \mathcal{C}^f(]0,T],\R^{n+1})$ such that $\alpha_{m} \longrightarrow \alpha$ for the Skorokhod topology. 
\\
Then for all $ A \text{ compact set in } ]0,T_{c}]$, $ \parallel \alpha_{m} -\alpha \parallel_{A}\longrightarrow 0 $, where $\parallel f\parallel_{A} = \sup_{t \in A} \parallel f(t) \parallel $.

Let $A$ be a compact set in $]0,T_{c}]$, then there exists a Lipschitz constant $ C_{A}$ of $F^{-1}$ in $\bigsqcup_{t \in A }M_{t} $, such that for all t in $A$,
$$d_{g(o)} (F^{-1}(\alpha_{m}(t)),F^{-1}(\alpha(t))) \le C_{A} \parallel \alpha_{m}(t)-\alpha(t) \parallel,$$
where $d_{g(o)}(x,y)$ is the distance in $M$ beetwen $x$ and $y$ for the metric $g(0)$. We also define $ d_{g(o),A} (f,g)= \sup_{t \in A} d_{g(o)}(f(t),g(t))$, where $f,g$ are $M$-valued function. We get:
$$d_{g(o),A} (\Psi(\alpha_{m}),\Psi(\alpha) \le C_{A} \parallel \alpha_{m}-\alpha \parallel_{A}.$$
So $\Psi(\alpha_{m}) \longrightarrow \Psi(\alpha)$ uniformly in all compact, so for the Skorokhod topology in $\mathcal{C}(]0,T_{c}],M )$.
\end{proof}

Let:
$$\tilde{Y}^\epsilon_{t} = (Y^\epsilon_{t} - Y^\epsilon_{0}) + (Y^\epsilon_{0} \1_{[\epsilon,T_{c}]}(t) +  \1_{[0,\epsilon]}(t) F(T_{c}-t,x_{o})).$$
Proposition \ref{tight_Y} gives the tightness of $Y^\epsilon_{t} - Y^\epsilon_{0}$, and $Y^\epsilon_{0} \1_{[\epsilon,T_{c}]}(t) +  \1_{[0,\epsilon]}(t) F(T_{c}-t,x_{o}) $ is a non-random sequence of functions that converges uniformly, hence $\tilde{Y}^\epsilon$ is  tight.
For strictly positive time $t$,
 $$ X^\epsilon_t= F^{-1}(T_{c}-t,\tilde{Y}^\epsilon_t) .$$
The previous lemma \ref{lem_sko} yields the tightness of $ X^\epsilon$.
Hence we have shown that:
 $$\forall \phi = (\epsilon_k)_k  \rightarrow 0 , \quad \exists  X^\phi_{]0,T_{c}]} , \qquad  X^{\epsilon_k}_{]0,T_{c}]} \overset{\mathcal{L}}{\rightarrow} X^\phi_{]0,T_{c}]} \text{ for an extracted sequence.} $$ 

\begin{proposition} \label{existance}
 Let $ \phi = (\epsilon_k)_k  \rightarrow 0 , \text{ and  } X^\phi_{]0,T_c]} ,s.t.  \quad  X^{\epsilon_k}_{]0,T_c]} \overset{\mathcal{L}}{\rightarrow} X^\phi_{]0,T_c]}$. Then  $ X^\phi_{]0,T_c]}$  is a  $\frac12 g(T_c-t)\text{-BM}$ in the following sense:
$$\forall \epsilon > 0 \quad X^\phi_{[\epsilon , T_c]}\overset{\mathcal{L}}{=}BM(\epsilon , X^\phi_{\epsilon}) $$.
\end{proposition}
\begin{proof}
 Let $\epsilon > 0$ then for large $k$:

 $$
\left\{ 
\begin{array}{l}
  X^{\epsilon_k} \text{  is a  } BM(\epsilon,X_\epsilon^{\epsilon_k})  \text{ after time } \epsilon \text{ , by Markov property} \\
                                                     \\
 \text{ and let  } X \text{ be a } BM(\epsilon,X^\phi_{_\epsilon})  \text{ after time } \epsilon 
\end{array}
\right.
$$
We want to show that $X = X^\phi $ after $\epsilon$ . 
So for sketch of the proof:
$$\begin{array}{c}
    X^{\epsilon_k} \overset{\mathcal{L}}{\underset{k \to \infty}{\longrightarrow}}X^\phi \\
\\
\text{so } X^{\epsilon_k}_{\epsilon} \overset{\mathcal{L}}{\underset{k \to \infty}{\longrightarrow}}X^\phi_{\epsilon} ,\\
\end{array}$$
we use the Skorokhod theorem, to have a $L_{2}$-convergence in a larger probability space: 
$$X^{'\epsilon_k}_{\epsilon} \overset{L_{2}, a.s.}{\underset{k \to \infty}{\longrightarrow}}X^{'\phi}_{\epsilon},$$
  with $X^{'\epsilon_k}_{\epsilon} \overset{\mathcal{L}}{=} X^{\epsilon_k}_{\epsilon} $
 and $ X^{'\phi}_{\epsilon} \overset{\mathcal{L}}{=} X^\phi_{\epsilon} $.
We use  convergence of solution of S.D.E with initial conditions converging in $ L_{2}$ (e.g. in  Stroock and Varadhan \cite{Stroock}), to get:
$$ BM(\epsilon,X^{'\epsilon_k}_{\epsilon})  \overset{\mathcal{L}}{\underset{k \to \infty}{\longrightarrow}} BM (\epsilon,X^{'\phi}_{\epsilon}),$$
 $$BM(\epsilon,X^{'\epsilon_k}_{\epsilon}) \overset{\mathcal{L}}{= }X^{\epsilon_k}_{ [\epsilon,T_c]} ,$$
 $$ BM(X^{'\phi}_{\epsilon}) \overset{\mathcal{L}}{= }BM(\epsilon,X^{\phi}_{\epsilon}).$$
We use that
$$ X^{\epsilon_k} \overset{\mathcal{L}}{\underset{k \to \infty}{\longrightarrow}} X^\phi $$
to conclude, after identification of the limit:
$$  X = BM(\epsilon,X^{\phi}_{\epsilon}) \overset{\mathcal{L}}{=}  X^\phi_{[\epsilon ,T_c]} .\\ $$
Hence the process  $ X^\phi$  is a  $\frac12 g(T_c-u)_{u \in ]0 ,T_c]}\text{-BM} $ in the above sense, we call it  "without birth" . 

\end{proof}

We now show that, in the sphere case, the $\frac12 g(T_c-u)_{u \in ]0 ,T_c]}\text{-BM} $ is, after a change of time, nothing else than a $BM(g(0))_{]-\infty , 0]}$, this will give uniqueness in law of such process.

\begin{proposition}\label{intro-mcf-sphere}
 Let $g(t)$ be a family of metrics which comes from a mean curvature flow on the sphere. Then the $ \tilde{g}(u)=\frac12 g(T_{c}-u)_{u \in ]0 ,T_{c}]}\text{-BM}$ is unique in law.
\end{proposition}
\begin{proof}
Let $R_{0}$ be the radius of the first sphere. Then $T_{c}= \frac{R_0^2}{2n}$, and by direct computation we obtain:
$$F(t,x) = \frac{\sqrt{R_0^2-2nt}}{R_0} x ,$$
$$g(t) = \frac{R_0^2-2nt}{R_0^2} g(0).$$
So for all $ f \in \mathcal{C}^{\infty}(S) $ we have:
$$ \D_{g(t)}f = \frac{R_0^2}{R_0^2-2nt} \D_{g(0)}f $$
and
$$
\begin{array}{rcl}
\nabla^{g(s)}df(X_{i},X_{j}) &=& f_{ij}-\Gamma_{ij}^{k}(s,.)f_{k} \\
                        &=& f_{ij}-\Gamma_{ij}^{k}(0,.)f_{k}  \text{ because the metrics are homothetic} \\
			&=& \nabla^{g(0)}df(X_{i},X_{j}).
\end{array}$$
Let $X $ be a $\frac12 g(T_{c}-u)_{u \in ]0 ,T_{c}]}\text{-BM} $.
For all $ f \in \mathcal{C}^{\infty}(S)$, $\ u \in ]0,T_{c}] $ and for all $T_{c}>t\ge u$ we have:

$\begin{array}{rcl}
f(X_{t})-f(X_{u}) &\overset{ \cal M }{\equiv}&  \frac12 \int_{u}^{t} \D _{\tilde{g}(s)}f(X_{s}) ds  \\
         &\overset{ \cal M }{\equiv}& \frac12 \int_{u}^{t} \nabla^{\tilde{g}(s)}df(*dX,*dX)  \\
	 &\overset{ \cal M }{\equiv}& \frac12 \int_{u}^{t} \nabla^{g(0)}df(*dX,*dX)  \\
 \end{array}$

\
\newline 
 $$df(X)_{]0,T_{c}]} \overset{d \cal M }{\equiv} \frac12  \nabla^{g(0)}df(*dX,*dX),$$
hence $X_{]0,T_{c}]}$ is a $g(0)$- martingale.
From \cite{metric}:
\begin{equation}\label{eq_sphere}
 df(X_{t}(x)) = \langle \nabla^{\tilde{g}(t)} f ,\tpar_{0,t}v_{i} \rangle_{\tilde{g}(t)}dW^{i} + \frac12 \D_{\tilde{g}(t)}(f)(X_{t}(x)) dt ,
\end{equation}

 with abusive notation (because we have no starting point, to get sense we have to take the conditional expectation at a time before $t$).
\\
It follows from (\ref{eq_sphere}):
$$ df(X_{t}(x)) =  \parallel \nabla^{\tilde{g}(t)} f (X_{t}(x))\parallel_{\tilde{g}(t)} dB_{t} + \frac12 (\frac{R_0^2}{R_0^2-2n(T_{c}-t)}) \D_{0}f(X_{t}(x)) dt, $$
where $B_{t}$ is some real-valued Brownian motion.
With help of the first computation,
$$ df(X_{t}) = \sqrt{\frac{R_0^2}{nt}} \parallel \nabla^{g(0)} f (X_{t})\parallel_{g(0)} dB_{t} + \frac12 (\frac{R_0^2}{nt}) \D_{0}f(X_{t}) dt.$$
Now consider the solution of:
$$\phi'(t)R_{0}^2 = n(\phi(t)) \text{  such that   } \phi(0)=T_{c}. $$
i.e. the function
 $$ \phi (t) =T_{c} e^{\frac{t}{2T_{c}}}.$$

 We get that $ X_{\phi (t)} = (BM_{g(0)})_{ t } $. According to the usual characterization of a Brownian motion \cite{Emery}.

So by this deterministic change of time, and by the uniqueness in law of a 
$ (BM_{g(0)})_{]-\infty,0]}$ on the sphere, we get the uniqueness in law of a $\frac12 g(T_{c}-u)_{u \in ]0 ,T_{c}]}\text{-BM} $ on a sphere.
\end{proof}

We have essentially used the conformality of this family of metric, that does not change the martingale family. Even if the beginning manifold is strictly convex, this is not the case in general. But we will see, in the next section, that the result is also true.

\section{Kendall-Cranston Coupling}
In this section the manifold $M$ is compact and strictly convex. The goal in this section is to prove the uniqueness in law of the $g(T_{c}-t)$-BM. This section will be cut in two parts, the first will be a geometric result inspired by the work of Huisken, the second will be an adaptation of the Kendall-Cranston coupling. We will, by a deterministic change of time, transform a $g(T_{c}-t)$-BM (the existence of which comes from proposition \ref{existance}) into a $\tilde{g}(t)_{]-\infty,0]}$-BM which has good geometric properties.

\begin{remark}
  In the two last sections in \cite{Huisken}, Huisken considers, like Hamilton for the Ricci flow, the normalized mean curvature flow. That consists in dilating the manifolds $M_{t}$ by a coefficient to obtain constant volume manifolds. He obtains a positive coefficient of dilation $\psi (t)$ that satisfies the following property.
\end{remark}

\begin{stheorem} \cite{Huisken} \label{H-Theorem}

For all $ t \in [0, T_{c}[ $, define $ \tilde{F}(.,t) =  \psi(t) F( . ,t)  $ such that $\int_{\tilde{M}_{t}} d\tilde{\mu}_{t} = |M_{0}|$, and $ \tilde{t} (t) = \int_{0}^{t} \psi^{2}(\tau) d\tau $, then there exist several positive constants $\delta, C$ such that:
\begin{enumerate}[i)]
  \item  $\tilde{T_{c}} = \infty  $
  \item   $\tilde{H}_{max}(\tilde{t})-\tilde{H}_{min}(\tilde{t}) \leq  Ce^{-\delta \tilde{t} } $
  \item  $|\frac{\partial}{\partial \tilde{t} } \tilde{g}_{ij}(\tilde{t})| \leq Ce^{-\delta \tilde{t} } $
   \item  $\tilde{g}_{ij}(\tilde{t}) \rightarrow \tilde{g}_{ij}(\infty)$ when  $\tilde{t}\rightarrow \infty $ uniformly, for the  $C^{\infty}-topology$, and the convergence is exponentially fast. 
   \item   $\tilde{g}(\infty)$ is a metric such that $(M,\tilde{g}(\infty))$ is a sphere.
\end{enumerate}
\end{stheorem}

We will now give the change of time propositions.

\begin{proposition}

Let $ \psi: [0,T_{c}[ \longrightarrow ]0,\infty[$ as above, $\tilde{t}$ defined by:
$$\begin{array}{rcl}
  \tilde{t}: [0,T_{c}[ &\longrightarrow& [0,\infty[ \\
   t &\longmapsto& \int_{0}^{t} \psi^{2}(\tau) d\tau, \\
 \end{array}$$
for all $ t \in [0 ,\infty[$, define
 $$ \quad \tilde{g}(t) = \psi^{2}(\tilde{t}^{-1}(t))g(\tilde{t}^{-1}(t)),$$
where $g(t)$ is the family of metrics coming from a mean curvature flow, and $ X_{t}$ is a $g(t)$-BM .
Then:
 
 $$ t \longmapsto X_{\tilde{t}^{-1}(t)}  \text{ is a } \tilde{g}(t)\text{-BM}  \text{ defined on }  [0,\infty[.$$
\end{proposition} 
\begin{proof}

Let $f \in \mathcal{C}^{\infty} (M) $:
 
 $$\begin{array}{rcl}
 f(X_{\tilde{t}^{-1}(t)}) & \overset{ \cal{M} }{\equiv} & \frac 12 \int_{0}^{\tilde{t}^{-1}(t)} \Delta_{g(s)}f(X_{s}) ds \\
   & \overset{ \cal{M} }{\equiv} & \frac 12 \int_{0}^{t} \Delta_{g(\tilde{t}^{-1}(s))}f(X_{\tilde{t}^{-1}(s)})  (\tilde{t}^{-1})'(s)ds \\
   &\overset{ \cal{M} }{\equiv} & \frac 12 \int_{0}^{t} \Delta_{\frac{1}{(\tilde{t}^{-1})'(s)} g(\tilde{t}^{-1}(s)) }f(X_{\tilde{t}^{-1}(s)})  ds .\\
   \end{array}$$ 
Using  
$$ \psi^{2}( \tilde{t}^{-1}(s)) (\tilde{t}^{-1})'(s) = 1 ,$$
we obtain:
$$\frac{1}{(\tilde{t}^{-1})'(s)} g(\tilde{t}^{-1}(s)) = \tilde{g}(s) .$$ 
\end{proof}

\begin{proposition}\label{change-time}

Let $X^{T_{c}}_{t} $, with $t \in ]0,T_{c}]$, be a  $g(T_{c}-t)\text{-BM}$. Let  $\tau$ be defined by:
$$\begin{array}{rcl}
  \tau: ]0,T_{c}] &\longrightarrow& ]-\infty,0] \\
   t  &\longmapsto& -\tilde{t}(T-t). \\
 \end{array}$$
Let $\tilde{g}(t)$ be defined by: 
$$\tilde{g}(t) = \psi^{2}(T_{c}-\tau^{-1}(t))g(T_{c}-\tau^{-1}(t)) \text{  } \forall t \in ]-\infty,0]. $$ 
Then: 

$$t \mapsto X^{T_{c}}_{\tau^{-1}(t)}  \text{ is a } \tilde{g}(t)\text{-BM}. $$ 

\end{proposition}

\begin{proof}

Let $f \in \mathcal{C}^{\infty} (M) $ and  $ s < t $,
$$\begin{array}{rcl}
 f(X^{T_{c}}_{\tau^{-1}(t)}) -  f(X^{T_{c}}_{\tau^{-1}(s)}) & \overset{ \cal{M} }{\equiv} & \frac 12 \int_{\tau^{-1}(s)}^{\tau^{-1}(t)} \Delta_{g(T_{c}-u)}f(X^{T_{c}}_{u}) du \\
   & \overset{ \cal{M} }{\equiv} & \frac 12\int_{s}^{t} \Delta_{g(T_{c}-\tau^{-1}(u))}f(X^{T_{c}}_{\tau^{-1}(u)})  (\tau^{-1}(u))'(s)du \\
   &\overset{ \cal{M} }{\equiv} & \frac 12 \int_{s}^{t} \Delta_{\frac{1}{(\tau^{-1})'(u)} g(T_{c}-\tau^{-1}(u)) }f(X^{T_{c}}_{\tau^{-1}(u)})  du .\\
\end{array}$$
We have $ -\tilde{t} (T_{c} -\tau^{-1}(u))= u $, and  
$$ (\tau^{-1})'(u) \psi^{2}(T_{c}-\tau^{-1}(u)) = 1. $$
We obtain 
 $$ f(X^{T_{c}}_{\tau^{-1}(t)}) -  f(X^{T_{c}}_{\tau^{-1}(s)})  \overset{ \cal{M} }{\equiv}  \frac 12 \int_{s}^{t} \Delta_{ \psi^{2}(T_{c}-\tau^{-1}(u)) g(T_{c}-\tau^{-1}(u)) }f(X^{T_{c}}_{\tau^{-1}(u)})  du $$
 
i.e.
 $$f(X^{T_{c}}_{\tau^{-1}(t)}) -  f(X^{T_{c}}_{\tau^{-1}(s)})  \overset{ \cal{M} }{\equiv}  \frac 12 \int_{s}^{t} \Delta_{ \tilde{g}(u) }f(X^{T_{c}}_{\tau^{-1}(u)})  du . $$

\end{proof}
\begin{remark}
 By the above theorem \ref{H-Theorem}, we know that $\tilde{g}(t)$ tends to a sphere metric as $t$ goes to $-\infty$. The above proposition transforms  ``two'' $g(T_{c}-t)$-BM into ``two'' $\tilde{g}$-BM so we will use the standardization of the metric into sphere metric and also the large time interval to perform the coupling.
\end{remark}

Let $\tau_{x}$ be a plane in $T_{x}M$ and $g(t)$ be a metric over $M$,  we denote by $K (t,\tau_{x})$ the sectional curvature of the plane $\tau_{x}$ according to the metric $g(t)$.
We will now give a few geometric lemmas that will be used later, for simplicity we will take positive times.
\begin{lemma}\label{geom-lemme}
\label{lem1}

  Let $g(t)$ be a family of metrics on a manifold $M$, and  $g(\infty)$ a metric that makes $M$ into a sphere, suppose that:
  \begin{enumerate}[i)]
 \item  $g(t) \longrightarrow g(\infty) \text{ uniformly, when } t \longrightarrow \infty $
   $   \text{ for the  } C^{\infty}-topology $ exponentially fast, i.e.: 
   $  \forall n \in \N , \forall \text{ multi-indices } ({i_{1}},...,{i_{k}}) \text{ such that } \sum i_{k}=n,\exists C_{n},\delta_{n} >0  , \text{ such that: } $
   $$ |\frac{\partial^{n}}{\partial X_{i_{1}}..X_{i_{k}} }g_{ij}(t) - \frac{\partial^{n}}{\partial X_{i_{1}}..X_{i_{k}} } g_{ij} (\infty)| \leq C_{n}e^{-\delta_{n} t }$$ 
 \item $\exists \delta, C^{1} > 0 $ such that  $|\frac{\partial}{\partial t } g_{ij}(t)| \leq C^{1}e^{-\delta t } $
 \item  $ vol_{g(t)}(M) = vol_{g(0)}(M)$
   \end{enumerate}
Then:

 for all $ \epsilon >0 $ , there exists $  T \in [0,\infty [ $, $\exists C ,cst, cst_{1}\in \R^{+} $ and $c_{n}(cst,V) > 0$ such that, $ \forall t \in [T , \infty [ $ the following conditions are satisfied:
 \begin{enumerate}[i)]  
   \item
     for all $x$ in $M$ and  for all plane $\tau_{x} \subset T_{x}M$,\quad $\mid K (t,\tau_{x})- cst \mid \le \epsilon$.
   \item 
     $ \vert \rho_{t}-\rho_{\infty} \vert_{M\times M  } \leq cst_{1}e^{-\delta t}$.
   \item
      $ \rho'_{t}(x,y):= \frac{d}{dt} \rho_{t}(x,y)\leq C \text{  in a compact CC of } M\times M , $ 
 
\end{enumerate}
  where the constant $ cst $, comes from the radius of $M$ with respect to $g(\infty)$, $\rho_{t}(x,y)$ is the distance between $x$ and $y$ for the metric $g(t)$, and 
$$CC  =\{(x,y)\in M \times M, s.t. \quad \rho_{t}(x,y) \leq min(\frac{\pi}{2\sqrt{(cst + \epsilon)}}, c_{n}(cst,V)) , \forall t > T \}.$$
\end{lemma}
\begin{proof}

 Let us prove i).
\\
Curvatures are functions of second order derivatives of the metric tensor.
We give the definitions of curvatures tensors, to make this point clear. Conventions are as in  \cite{Lee},\cite{Rjost},\cite{Jost}, in particular, we use Einstein's summation convention..
For a metric connection without torsion (Levi-Civita connection), we recall standard definitions:
\\
-the Christoffel symbols:
  $$\Gamma _{ij}^{k} = \frac 12 g^{kl} ( \frac{\partial}{\partial x_{i}} g_{jl} + 
                                        \frac{\partial}{\partial x_{j}} g_{il} -
  		                        \frac{\partial}{\partial x_{l}} g_{ij}  )$$
-the  (3,1) Riemann tensor:
   $$ R(X,Y)Z = \nabla_{X}\nabla_{Y} Z -\nabla_{Y}\nabla_{X} Z - \nabla_{[X,Y]} Z $$
-the (4,0) curvature tensor: 
   $$  R_{m}(X,Y,Z,W) = \langle R(X,Y)Z , W \rangle $$
-the sectional curvature:
   $$ K(X,Y) = \frac{R_{m}(X,Y,Y,X)}{|X|^{2}|Y|^{2}-\langle X,Y \rangle^{2} } $$
We see that the sectional curvature depends on the metric and its derivatives up to order two, so 
$ \forall x \in M,$ for all plane $\tau_{x} \subset T_{x}M$,  
$$ \lim_{t \rightarrow \infty} K(t,\tau_{x}) = cst .$$
Also, for all $ \epsilon>0$, there exists $  T $ such that  $\forall t > T $, for all $x$ in $M$ and for all plane $\tau_{x} \subset T_{x}M$,
 $$\mid K (t,\tau_{x})- cst \mid \le \epsilon.$$

For the third point iii):
\\
for $(x,y) \in CC$, where $CC$ is defined above, we will show that we have the uniqueness of minimal $g(t)$-geodesic from $x$ to $y$, for all time $t>T$, because we have the well-known Klingenberg's result (e.g. \cite{GHL} page 158) about injectivity radius of compact manifold whose sectional curvature is bounded above. To use Klingenberg's lemma, we have to bound the shortest length of a closed geodesic. We will use Cheeger's theorem  page 96 \cite{cheeger}. Since by the convergence of the metric, we have the convergence of the Ricci curvature, we obtain that they are bounded by the same constant. We obtain, using Myers' theorem that all diameters are then bounded above. The volumes are constant so bounded below, all sectional curvatures of $M$ are bounded in absolute value from above. So by Cheeger's theorem there exists a constant $c_{n}(K,d,V) > 0$ that bounds the length of smooth closed geodesics. Hence, for large time, using Klingenberg's lemma, we get a uniform bound , in time, of the injectivity radius (i.e $ min(\frac{\pi}{2\sqrt{(cst + \epsilon)}}, c_{n}(cst,V)) $).

So for all $t>T$, there exist only one $g(t)$-geodesic between $x$ and $y$, we denote it  $\gamma^{t} $. Let $E(\gamma^{t}) =  \int_{0}^{1} \langle  \dot{\gamma}^{t}(s) , \dot{\gamma}^{t}(s) \rangle_{g(t)} ds$ be the energy of the geodesic where $\dot{\gamma}^{t}(s)= \frac{\partial}{\partial s} \gamma^{t}(s)$, $\rho^{2}_{t}(x,y) = E(\gamma^{t})$. We compute:

$$\begin{array}{rcl}
 2(\frac{\partial}{\partial t}|_{t=t_{0}}\rho_{t}(x,y)) (\rho_{t}(x,y)) &=& \frac{\partial}{\partial t}|_{t=t_{0}} E(\gamma^{t}) \\
&=&\int_{0}^{1} \langle  \dot{\gamma}^{t_{0}}(s) , \dot{\gamma}^{t_{0}}(s) \rangle_{\frac{\partial}{\partial t}|_{t=t_{0}} g(t)}ds \\
 &+&   2  \int_{0}^{1} \langle  D_{t}|_{t=t_{0}}\frac{\partial}{\partial s}\gamma^{t}(s) , \frac{\partial}{\partial s}\gamma^{t_{0}}(s) \rangle_{g(t_{0})} ds                                     \\
 &=& \int_{0}^{1} \langle  \dot{\gamma}^{t_{0}}(s) , \dot{\gamma}^{t_{0}}(s) \rangle_{\frac{\partial}{\partial t}|_{t=t_{0}} g(t)}ds \\
 &+&   2  \int_{0}^{1} \langle D_{s} \frac{\partial}{\partial t}|_{t=t_{0}}\gamma^{t}(s) , \frac{\partial}{\partial s}\gamma^{t_{0}}(s) \rangle_{g(t_{0})} ds                                     \\
  \end{array}$$

Let $ X =  \frac{\partial}{\partial t}|_{t=t_{0}} \gamma^{t}(s)$  be a vector field such that $ X(x) = 0_{T_{x}M} ,X(y) = 0 _{T_{y}M}$, because we do not change the beginning and terminal point.
The covariant derivative is computed with the Levi-Civita connection associated to $g(t_{0})$. Hence we obtain:
$$ \int_{0}^{1} \langle D_{s} \frac{\partial}{\partial t}|_{t=t_{0}}\gamma^{t}(s) , \frac{\partial}{\partial s}\gamma^{t_{0}}(s) \rangle_{g(t_{0})} ds                                   = \int_{0}^{1} \langle \nabla_{\dot{\gamma}^{t_{0}}(s)} X, \frac{\partial}{\partial s}\gamma^{t_{0}}(s) \rangle_{g(t_{0})} ds,$$
also:
 $$\langle \nabla_{\dot{\gamma}^{t_{0}}(s)} X, \frac{\partial}{\partial s}\gamma^{t_{0}}(s) \rangle_{g(t_{0})} = \frac{\partial}{\partial s}\langle  X, \frac{\partial}{\partial s}\gamma^{t_{0}}(s) \rangle_{g(t_{0})}, $$
because the connection is metric and $\gamma^{t_{0}}$ is a $g(t_{0})$-geodesic. Hence 
$$\int_{0}^{1} \frac{\partial}{\partial s}\langle  X, \frac{\partial}{\partial s}\gamma^{t_{0}}(s) \rangle_{g(t_{0})} ds=[\langle  X, \frac{\partial}{\partial s}\gamma^{t_{0}}(s) \rangle_{g(t_{0})} ]_{0}^{1}=0.$$
Finally, we obtain:
 \begin{equation} \label{eq1}
  \frac{\partial}{\partial t}|_{t=t_{0}}\rho_{t}(x,y) = \frac{1}{2 \rho_{t_{0}}(x,y) } \int_{0}^{1} \langle  \dot{\gamma}^{t_{0}}(s) , \dot{\gamma}^{t_{0}}(s) \rangle_{\frac{\partial}{\partial t}|_{t=t_{0}} g(t)}ds  .
 \end{equation}
We will now control the second term in the previous equation.
By the exponential convergence of the metric, we could assume that the time is in the compact interval $[0,1]$. The manifold is compact, so we have a finite family of charts (indeed, we may assume that we have two charts, because the manifold has a metric which turns it into a sphere). The support of this chart could be taken to be relatively compact, and in this chart we can take the Euclidien metric i.e $\langle  \partial_{i} ,\partial_{j} \rangle_{E} = \delta_{i}^{j}$. This is not in general a metric on $M$. For the simplicity of expression, after taking the minimum over all charts we may assume that we just have one chart.
Let $ S_{1}$ be a sphere in $\R^{n} $ with the Euclidean metric. The functional:
  $$\begin{array}{rcl}
  [0,1] \times S_{1} \times M &\longrightarrow& \R \\
  (t,v,x) &\longmapsto& g_{ij}(t , x)v_{i}v_{j}\\
   \end{array}$$
  reaches its minimum $C>0$, so: 
  $$\Vert T \Vert _{E} \leq C^{-1}\Vert T \Vert_{g(t)} , \forall t \in [0,1] , \forall T \in TM. $$ 
Hence, for the equation (\ref{eq1}) we get the estimate:
 $$\begin{array}{rcl}
 \Big{|}\frac{\partial}{\partial t}|_{t=t_{0}}\rho_{t}(x,y)\Big{|} &\leq& \frac{1}{2 \rho_{t_{0}}(x,y) } C^{1} e^{- \delta t_{0}}\int_{0}^{1} \Big{|} \langle \dot{\gamma}^{t_{0}}(s) , \dot{\gamma}^{t_{0}}(s) \rangle_{E}\Big{|}ds \\
 & \leq& \frac{1}{2 \rho_{t_{0}}(x,y) } C^{1} (C)^{-1} e^{- \delta t_{0}} \int_{0}^{1} \Big{|} \langle \dot{\gamma}^{t_{0}}(s) , \dot{\gamma}^{t_{0}}(s) \rangle_{g(t_{0})}\Big{|}ds \\
 & \leq& \frac{1}{2 } C^{1} (C)^{-1} e^{- \delta t_{0}}. \\
\end{array}$$
This expression is clearly bounded.

For the second point ii),
\\
 let $x,y \in M $ take  $\gamma_{\infty}$ be a $g(\infty) $-geodesic that joins $x$ to $y$. Then we have, on the one hand,
$$\begin{array}{rcl}
   \rho^{2}_{t}(x,y)-\rho^{2}_{\infty}(x,y) &\leq& \int_{0}^{1} \langle \dot{\gamma}_{\infty}(s) , \dot{\gamma}_{\infty}(s) \rangle_{g(t) - g(\infty)} ds \\
   
    &\leq& Cst e^{-\delta t} \int_{0}^{1} \Vert \dot{\gamma}_{\infty}(s) \Vert^{2} _{g(\infty)} ds \\
    &\leq& Cst e^{-\delta t} diam^{2}_{g(\infty)}(M);\\
    \end{array}$$
where the constant changes and depends on the previous constant.
On the other hand, we have:
   
   $$\begin{array}{rcl}
   \rho_{\infty}^{2}(x,y)-\rho^{2}_{t}(x,y) &\leq& \int_{0}^{1} \langle \dot{\gamma}^{t}(s) , \dot{\gamma}^{t}(s) \rangle_{g(\infty)-g(t)} ds \\
    &\leq& Cst e^{-\delta t} \int_{0}^{1} \Vert \dot{\gamma}^{t}(s) \Vert^{2} _{g(t)} ds \\
    &\leq& Cst e^{-\delta t} diam^{2}_{g(t)}(M)\\
    &\leq& cst_{1} e^{-\delta t} ,\\
    \end{array}$$
for some constant $cst_{1}$ , and we use Myers theorem for the last inequality to get a uniform upper bound of the diameter (because all Ricci curvature are uniformly bounded). 
\\
We get exponential convergence of the length.

\end{proof}


We will now show uniqueness in law of a $g(T_{c}-t)$-BM. By proposition \ref{change-time}, this uniqueness is equivalent to uniqueness in law of a $\tilde{g}(t)_{]-\infty,0]}$-BM. This family of metrics, $\tilde{g}(t)$, satisfies:
  $$\tilde{g}(t) \longrightarrow \tilde{g}(-\infty) \text{  for the  }C^{\infty}\text{-topology} .$$
Let $Z^{1}, Z^{2}$ be two $\tilde{g}\text{-BM}_{]-\infty , 0]} $ and  $N << T$ where $T$ is the time of the lemma \ref{geom-lemme}, i.e the time up to which all bounds of the lemma are under control. Geometry before this time is similar to the geometry of the sphere. 
So the result of uniqueness in law for Brownian motion defined in a product probability space, indexed by $\R$ in a compact manifold (e.g. \cite{Emery_Schachermayer},\cite{Marcfiltration}) could give the heuristics to our results.
As we can see in \cite{metric} the $g(t)$-stochastic development and the $g(t)$-horizontal lift of a $g(t)$-BM is well defined.
\\
We will consider a new process  $ Z^{3}_{N,t}$ equal in law to $  Z^{2}$ after $N$ and equal to $  Z^{2}$ before. In the sequel we will note  
$Z^{3}_{t}$ for $Z^{3}_{N,t}$. The construction, after time $N$, will be given by localization in a stochastic interval.

Let $ T^{N}_{0} = N, $ and for all $ t \leq N $, $Z^{3}_{N,t} = Z^{2}_{t} $.\\
 1)  we will let $Z^{3}_{t}$ evolve independently of $Z^{1}_{t}$ i.e.
 $Z^{3}_{t}$ is a $g(T^{N}_{0}+ .)$-BM which starts at $Z_{ T^{N}_{0} }^{3}$ and the $\R^{n}$-valued Brownian motion that drives $Z^{3}_{t}$ will be independent with the one that drives $Z^{1}_{t}$.

Let $T^{N}_{1} = (N+\frac12) \wedge inf \{ t> T^{N}_{0}, \text{  } \rho_{t}(Z^{1}_{t},Z^{3}_{t}) \leq  \frac{ \frac{\pi}{\sqrt{cst +\epsilon}} \wedge \frac{C_{n}(d,K,cst-\epsilon) }{2}} {4}\} \wedge T $.
The constant $\epsilon$ is just taken to be small enough.

Let $C_{N} = inf \{ t> N, \text{ } Z^{1}_{t}=Z^{3}_{t}   \}. $
\\
2) At time $T^{N}_{1}$:
\begin{itemize}
\item if $\rho_{T^{N}_{1}}(Z^{1}_{T^{N}_{1}},Z^{3}_{T^{N}_{1}}) \leq  \frac{ \frac{\pi}{\sqrt{cst +\epsilon}} \wedge \frac{C_{n}(d,K,cst-\epsilon) }{2}} {4}$, these two points ($Z^{3}_{T^{N}_{1}}$ and $Z^{1}_{T^{N}_{1}}$) are close enough to make mirror coupling. The distance between these two points is strictly less than the injectivity radius $i_{g(t)}(M)$, hence we have uniqueness of the geodesic that joins these two points.
 After $T^{N}_{1} $ and before $C_{N}$, we build $Z^{3}_{t}$ as the $g(T^{N}_{1} + .)$-BM that starts at $ Z_{T^{N}_{1}}^{3} $, and solves:

 $$*dZ_{t}^{3} =U^{3}_{t}*d((U^{3}_{t})^{-1} m^{t}_{Z^{1}_{t},Z^{3}_{t}} U^{1}_{t} e_{i}dW^{i}_{t}) $$
and after $C_{N}$,
  $$ Z^{3}_{t} = Z^{1}_{t},\quad  C_{N} \leq t , $$

where  $U^{3}_{t}$ is the horizontal lift of $Z^{3}_{t}$, to be correct we have to express a system of stochastic differential equations as in Kendall \cite{Kendall}, $U^{1}_{t}$ is the horizontal lift of $Z^{1}_{t}$, and $dW^{i}_{t}$ are Brownian motion that drives $Z^{1}_{t}$, the mirror map $m^{t}_{x,y}$ consists in transporting a vector along the unique minimal $g(t)$-geodesic that joins $x$ to $y$ and then reflecting it in the hyperplane of $(T_{y}M, g(t))$ which is perpendicular to the incoming geodesic.

By isometry property of the horizontal lift of the $g(t)$-BM (see \cite{metric}), 
$$ (U^{3}_{t})^{-1} m^{t}_{Z^{1}_{t},Z^{3}_{t}} U^{1}_{t} dW^{i}_{t}, $$ is an $\R^{n}$-valued Brownian motion.

Let  $T^{N}_{2} = (T^{N}_{1} +\frac12) \wedge inf \{ t> T^{N}_{1}, \text{ } \rho_{t}(Z^{1}_{t},Z^{3}_{t}) > \frac{ \frac{\pi}{\sqrt{cst +\epsilon}} \wedge \frac{C_{n}(d,K,cst-\epsilon)}{2}}{2}  \} \wedge T \wedge C_{N} $.

\item if $\rho_{T^{N}_{1}}(Z^{1}_{T^{N}_{1}},Z^{3}_{T^{N}_{1}}) >  \frac{ \frac{\pi}{\sqrt{cst +\epsilon}} \wedge \frac{C_{n}(d,K,cst-\epsilon) }{2}} {4}$
then $T^{N}_{2} = T^{N}_{1}$.
\end{itemize}
Iterate step 1 and 2 successively (changing $T^{N}_{0}$ by $T^{N}_{2}$ and $T^{N}_{1}$ by $T^{N}_{3}$ in step 1, changing $T^{N}_{1}$ by $T^{N}_{3}$ and $T^{N}_{2}$ by $T^{N}_{4}$ in step 2 ..., after time $T$ if we have no coupling, we let $Z^{3}$ evolve independently of $Z^{1}_{t}$ until the end), we build by induction the process $Z^{3}_{t} $ and a sequence of stopping times.
We sketch it as:
\begin{itemize}
 \item if $C_{N} < T$
$$T^{N}_{0} \stackrel{\text{ independent }}{ \longrightarrow} T^{N}_{1} \stackrel{coupling}{ \longrightarrow} T^{N}_{2} \stackrel{independent}{ \longrightarrow} T^{N}_{3} \stackrel{coupling}{ \longrightarrow} T^{N}_{4} \text{  ...  } C_{N} \stackrel{Z^{3}_{t} = Z^{1}_{t}}{ \longrightarrow} 0 $$

 \item if $C_{N} > T$
$$T^{N}_{0} \stackrel{\text{ independent }}{ \longrightarrow} T^{N}_{1} \stackrel{coupling}{ \longrightarrow} T^{N}_{2} \stackrel{independent}{ \longrightarrow} T^{N}_{3} \stackrel{coupling}{ \longrightarrow} T^{N}_{4} \text{  ...  } T \stackrel{independent}{ \longrightarrow} 0 $$

\end{itemize}
   
\begin{proposition}
 The two processes $ Z^{3}$ and $Z^{2}$ are equal in law.  
\end{proposition}

 \begin{proof}
 It is clear that before $N$ the two processes are equal so equal in law.
After: 

 $Z_{N}^{3} = Z_{N}^{2}$.

 $$ \left\lbrace  \begin{array}{l}
       *dZ_{t}^{3} = \sum_{i} U^{3}_{t}e_{i} *dB^{i} , \text{ when } t \in [T^{N}_{2k},T^{N}_{2k+1}  ] \text{ and }  T^{N}_{2k+1} \leq  C_{N}  \\
       
       *dZ_{t}^{3} =\sum_{i} U^{3}_{t}*d((U^{3}_{t})^{-1} m^{t}_{Z^{1}_{t},Z^{3}_{t}} U^{1}_{t}) e_{i}dW^{i}_{t}  , \text{ when } t \in [T^{N}_{2k+1},T^{N}_{2k+2}  ] ,\text{ and }  T^{N}_{2k+2} \leq  C_{N}\\
        
        Z^{3}_{t} = Z^{1}_{t} ,  C_{N} \leq t  \\
        \end{array} \right .$$
We write:
 $$\begin{array}{lcl}
 Z^{3}_{t} &= &  \sum_{k=0}^{\infty} \1 _{[T^{N}_{k},T^{N}_{k+1}]}(t)*d Z^{3}_{t}\\
             &=&  \sum_{k: even}.... + \sum_{k: odd} \\
 \end{array}$$
 \\
 Let $f \in C^{\infty}(M) $ then we have:
 \\
 for even $k$:
 $$ df(\1 _{[T^{N}_{k},T^{N}_{k+1}]}(t)*dZ^{3}_{t})  \overset{d \cal M }{\equiv} \frac12 \1 _{[T^{N}_{k},T^{N}_{k+1}]}(t)\Delta_{\tilde{g}(t)} f (Z^{3}_{t} ) dt $$
 for odd $k$:
 \begin{eqnarray*}   
 df(\1 _{[T^{N}_{k},T^{N}_{k+1}]}(t)*dZ^{3}_{t})   &=& \frac12 \1 _{[T^{N}_{k},T^{N}_{k+1}]}\Delta_{\tilde{g}(t)}f(Z^{3}_{t})dt\\		     
  \end{eqnarray*}
  
So  $ Z^{3}$ and $Z^{2}$ are two diffusions with the same starting distribution and the same generator, hence they are equal in law. For the gluing with $Z^{1}$ after $ C_{N} $ this is just the strong Markov property for $(t,Z)$.
 \end{proof}


 \begin{proposition}\label{alfa1}
 There exists $ \alpha >0$  such that:
 $$ \P ( T^{N}_{1} - N < \frac12 ) > \alpha  $$
 \end{proposition}
\begin{proof}
 By the $ C^{\infty}$-convergence of the metric we get:
$$\forall t < T ,  \vert \D_{\tilde{g}(t)} f -\D_{\tilde{g}(-\infty)} f\vert \leq \tilde{C}e^{\delta t} $$
where the constant comes from Theorem \ref{H-Theorem}, and the derivative of $f$ up to order two. We also obtain, by lemma \ref{geom-lemme}, for a constant $\epsilon_{2}$ that will be fixed below:
 $$ \vert \rho_{t}- \rho_{-\infty} \vert \leq \epsilon_{2} .$$
Over the sphere $ (M,\tilde{g}(-\infty))$, we have by ordinary comparison theorem:
$$\D_{\tilde{g}(-\infty)} \rho_{-\infty} (x) \leq (n) cot(\rho_{-\infty} (x)). $$
We can suppose after normalization that the radius of the sphere $ (M,\tilde{g}(-\infty))$ is one, $Radius_{-\infty}( M) = 1 $ (i.e. $cst = 1 $) in \ref{geom-lemme}. We deduce from above that:
 $$\D_{\tilde{g}(t)} \rho_{-\infty} (x) \leq (n) cot(\rho_{-\infty} (x)) + \tilde{C}e^{\delta t} .$$
In $ [N, T^{N}_{1}[$, we have $\rho_{t}(Z^{1}_{t},Z^{3}_{t}) > \frac{ \frac{\pi}{\sqrt{1 +\epsilon}} \wedge \frac{C_{n}(d,K,cst-\epsilon) }{2}} {4} $ so:
$$ \frac{ \frac{\pi}{\sqrt{1 +\epsilon}} \wedge \frac{C_{n}(d,K,cst-\epsilon) }{2}} {4}-\epsilon_{2} \leq \rho_{t}(Z^{1}_{t},Z^{3}_{t})-\epsilon_{2} \leq \rho_{-\infty}(Z^{1}_{t},Z^{3}_{t}) \leq \pi $$
We can choose $ \epsilon , \epsilon_{2}$ such that, $\frac{ \frac{\pi}{\sqrt{1 +\epsilon}} \wedge \frac{C_{n}(d,K,cst-\epsilon) }{2}} {4} -\epsilon_{2} \geq \beta > 0 $.
We obtain:
 $$ cot(\rho_{-\infty}(Z^{1}_{t},Z^{3}_{t})) \leq cot(\beta), $$
and 
 $$\D_{\tilde{g}(t)}\rho_{-\infty}(Z^{1}_{t},.)  ( Z^{3}_{t}) \leq (n) cot(\beta)+\tilde{C}e^{\delta T}, $$
(recall $T<<0$)
The progression of $ Z^{3}$ and $Z^{1}$ are independent between  $ [N, T^{N}_{1}]$ hence:
  $$ (Z^{1}_{t},Z^{3}_{t}) \text{ is a diffusion with generator  } \frac12( \D_{\tilde{g}(t),1}+\D_{\tilde{g}(t),2})$$
  i.e. 
  $$d\rho_{-\infty}(Z^{1}_{t},Z^{3}_{t}) = dM_{t} +\frac12( \D_{\tilde{g}(t)}\rho_{-\infty}(Z^{1}_{t},.)(Z^{3}_{t})+\D_{\tilde{g}(t)} \rho_{-\infty}(.,Z^{3}_{t})(Z^{1}_{t}))dt $$
  where $M_{t}$ is a local martingale, so
$$d\rho_{-\infty}(Z^{1}_{t},Z^{3}_{t}) \leq dM_{t} + (cot(\frac{\pi}{8}) + \tilde{C}e^{\delta T}) dt.  $$.

Let us compute the quadratic variation of this local martingale, i.e:
$$d\langle M, M \rangle_{t}= d\rho_{-\infty}(Z^{1}_{t},Z^{3}_{t}) d\rho_{-\infty}(Z^{1}_{t},Z^{3}_{t}),$$
with: 
\begin{equation} \label{dif_sto1}
  d\rho_{-\infty}(Z^{1}_{t},Z^{3}_{t})=d\rho_{-\infty}(Z^{1}_{t},.)*d Z^{3}_{t} + d\rho_{-\infty}(.,Z^{3}_{t})*d Z^{1}_{t}.
\end{equation}
Let $\gamma_{-\infty}(Z^{3}_{t},Z^{1}_{t}) (s)$ be the minimal $\tilde{g}(-\infty)$-geodesic beetwen  $ Z^{3}_{t}$ and $Z^{1}_{t} $ that exists and is unique almost everywhere because $Cut_{-\infty}(M) $ is a null measure subspace. We denote:
$$ v_{t}^{1}=\frac{\dot{\gamma}_{-\infty}(Z^{3}_{t},Z^{1}_{t}) (0)}{\Vert \dot{\gamma}_{-\infty}(Z^{3}_{t},Z^{1}_{t}) (0) \Vert_{\tilde{g}(-\infty)}} .$$ 
We complete $v_{t}^{1}$ with $v_{t}^{j}$  to get a  $\tilde{g}(-\infty)$-orthonormal basis. We rewrite $*dZ_{t}^{3}$ as:
 
\begin{eqnarray*}
*dZ_{t}^{3} &=& \sum U^{3}_{t}e_{i} *dB^{i} \\
            &=& \sum_{i,j} \langle U^{3}_{t}e_{i} , v_{t}^{j} \rangle_{\tilde{g}(-\infty)}v_{t}^{j}*dB^{i} \\
\end{eqnarray*}
Hence by Gauss lemma, we obtain:
\begin{eqnarray*}
d\rho_{-\infty}(Z^{1}_{t},.)*d Z^{3}_{t} &=& \sum d\rho_{-\infty}(Z^{1}_{t},.)U^{3}_{t}e_{i} *dB^{i}\\
           &=&  \sum_{i,j} d\rho_{-\infty}(Z^{1}_{t},.) \langle U^{3}_{t}e_{i} , v_{t}^{j} \rangle_{\tilde{g}(-\infty)}v_{t}^{j} *dB^{i} \\
	   &=& \sum_{i} d\rho_{-\infty}(Z^{1}_{t},.) \langle U^{3}_{t}e_{i} , v_{t}^{1} \rangle_{\tilde{g}(-\infty)}v_{t}^{1} *dB^{i} \\
	   &=&\sum_{i} \langle U^{3}_{t}e_{i} , v_{t}^{1} \rangle_{\tilde{g}(-\infty)} *dB^{i}. 
\end{eqnarray*}
It follows that:
$$ (d\rho_{-\infty}(Z^{1}_{t},.)*d Z^{3}_{t})(d\rho_{-\infty}(Z^{1}_{t},.)*d Z^{3}_{t}) = \sum_{i} \langle U^{3}_{t}e_{i} , v_{t}^{1} \rangle^{2}_{\tilde{g}(-\infty)} dt. $$
By the exponential convergence of the metric,
$$\langle U^{3}_{t}e_{i} , v_{t}^{1} \rangle_{\tilde{g}(-\infty)} \geq \langle U^{3}_{t}e_{i} , v_{t}^{1} \rangle_{\tilde{g}(t)} - \tilde{C}e^{\delta T}, $$
hence:
\begin{eqnarray*}
\sum_{i} \langle U_{t}e_{i} , v_{t}^{1} \rangle^{2}_{\tilde{g}(-\infty)} &\geq & 
\sum_{i} \langle U_{t}e_{i} , v_{t}^{1} \rangle^{2}_{\tilde{g}(t)} -2\tilde{C}e^{\delta T}
\sum_{i} \langle U_{t}e_{i} , v_{t}^{1} \rangle_{\tilde{g}(t)} + n (\tilde{C}e^{\delta T})^{2} \\
& = & \Vert  v_{t}^{1} \Vert_{\tilde{g}(t)}^{2} -2\tilde{C}e^{\delta T}
\sum_{i} \langle U_{t}e_{i} , v_{t}^{1} \rangle_{\tilde{g}(t)} + n (\tilde{C}e^{\delta T})^{2}\\
&\geq &\Vert v_{t}^{1}\Vert_{\tilde{g}(t)}^{2} -2\tilde{C}e^{\delta T}
n\Vert v_{t}^{1}\Vert_{\tilde{g}(t)} + n (\tilde{C}e^{\delta T})^{2} \text{  Schwartz}\\
&\geq & (\Vert v_{t}^{1}\Vert_{\tilde{g}(-\infty)}-\tilde{C}e^{\delta T})^{2} -2\tilde{C}e^{\delta T}n (\Vert v_{t}^{1}\Vert_{\tilde{g}(-\infty)}+\tilde{C}e^{\delta T} )\\
&+& n (\tilde{C}e^{\delta T})^{2} \\
&\geq & 1 - \tilde{C}e^{\delta T} (2-\tilde{C}e^{\delta T} +2(n+n\tilde{C}e^{\delta T} )-n\tilde{C}e^{\delta T} )\\
&\geq& \frac12 \text{  for a small enough T. }
\end{eqnarray*}
The independence of $Z^{1}_{t}$ et $Z^{3}_{t}$ gives,
\begin{eqnarray*}
 d\langle M_{t},M_{t} \rangle &=& (d\rho_{-\infty}(Z^{1}_{t},.)*d Z^{3}_{t})(d\rho_{-\infty}(Z^{1}_{t},.)*d Z^{3}_{t}) \\
& & +(d\rho_{-\infty}(.,Z^{3}_{t}))*d Z^{1}_{t})(d\rho_{-\infty}(.,Z^{3}_{t}))*d Z^{1}_{t}) 
\end{eqnarray*}
hence
$$ d\langle M_{t},M_{t} \rangle \geq 1 dt. $$
\\
\\
For simplicity write $\theta = \frac{ \frac{\pi}{\sqrt{1 +\epsilon}} \wedge \frac{C_{n}(d,K,cst-\epsilon) }{2}} {4}$, it follows from (\ref {dif_sto1}) that: 
\begin{align*}
& \P ( T^{N}_{1} - N < \frac12 )\\
 &=  \P (\exists t \in [N,N+1/2] \text{ s.t. } \rho_{t}(Z^{1}_{t},Z^{3}_{t})\leq  \theta \\
 & \geq  \P (\exists t \in [N,N+1/2] \text{ s.t. } \rho_{-\infty }(Z^{1}_{t},Z^{3}_{t}) \leq \theta -\epsilon_{2}) \\
 & \geq  \P (\exists t \in [N,N+1/2] \text{ s.t. } \pi+M_{t} + (cot(\beta) + \tilde{C}e^{\delta T})(t-N)\leq \theta-\epsilon_{2}) \\
 &\geq \alpha >0.\\
\end{align*}
For the last step, we use the usual comparison theorem for stochastic processes (e.g. Ikeda and Watanabe \cite{IW}).
\end{proof}
\\
We will now show that the coupling can occur between  $[T^{N}_{1},T^{N}_{2}] $ in a time smaller than $\frac12 $.\\

\begin{proposition}\label{alfa2}
There exists $ \tilde{\alpha} >0  $  such that:
 $$ \P ( C_{N}  < (T^{N}_{1} + \frac12 )\wedge T^{N}_{2} ) > \tilde{\alpha}.$$
\end{proposition}
\begin{proof}
 Between the two times $T^{N}_{1} $ and $T^{N}_{2} $, we have mirror coupling between $Z^{1}_{t}$ and $Z^{3}_{t}$. As in \cite{Kendall, Cranston} we have:
 $$d\rho_{t}(Z^{1}_{t},Z^{3}_{t}) = \rho'_{t}(Z^{1}_{t},Z^{3}_{t})dt + 2d\beta_{t}+\frac12 \sum_{i=2}^{n} I^{t}(J_{i}^{t},J_{i}^{t}) dt $$
$$*dZ_{t}^{3} =U^{3}_{t}*d((U^{3}_{t})^{-1} m^{t}_{Z^{1}_{t},Z^{3}_{t}} U^{1}_{t} e_{i}dW^{i}_{t}), $$
Where:\\
-$\beta_{t}$ is a standard real Brownian motion.
 \\
 -$\gamma_{t}(Z^{1}_{t},Z^{3}_{t}) (s)$ the minimal $\tilde{g}(t)$ geodesic between  $ Z^{1}_{t}$ and $Z^{3}_{t}. $
\\
 -$(\dot{\gamma}(Z^{1}_{t},Z^{3}_{t})(0), e_{i}(t))$ a $\tilde{g}(t)$-orthonormal basis of $T_{Z^{1}_{t}}M$.
 \\
 -  $J_{i}^{t}(s)$ the Jacobi field along $\gamma_{t}$ for the metric  $\tilde{g}(t) $, with initial condition $J_{i}^{t}(0)=e_{i}(t)$ and $J_{i}^{t}(\rho_{t}(Z^{1}_{t},Z^{3}_{t}))= \tpar_{\rho_{t}(Z^{1}_{t},Z^{3}_{t})}^{t,\gamma_{t}} e_{i}(t)$ i.e. the parallel transport for the metric $\tilde{g}(t) $ along $\gamma_{t}$, that is an orthogonal Jacobi field .
 \\
 -$I^{t}$ is the index bilinear form for the metric  $\tilde{g}(t). $
 \\
Between the times $T^{N}_{1} $ and $T^{N}_{2}$, we have:
 $$ \rho_{t}(Z^{1}_{t},Z^{3}_{t}) \leq  \frac{ \frac{\pi}{\sqrt{cst +\epsilon}} \wedge \frac{C_{n}(d,K,cst-\epsilon) }{2}} {2}$$
So by \ref{geom-lemme}, there exists a constant $C$ such that:
 $$ \rho'_{t}(x,y) \leq C.  $$
We have to show that between the times $T^{N}_{1} $ and $T^{N}_{2}$,
 $$\sum_{i=2}^{n} I^{t}(J_{i}^{t},J_{i}^{t})  $$ 
is bounded above.
We note   $r= \rho_{t}(Z^{1}_{t},Z^{3}_{t}) $, and $\gamma $ for $ \gamma^{t}$.
Let $ G(s)$ be a real-valued function and  $K_{i}^{t}$ be the orthogonal vector field over $\gamma $ defined by:
 $$K_{i}^{t}(s) = G(s) (\tpar_{t}^{\gamma_{t}} e_{i}(t))(s) $$
 where  $ G(0)= G(r) = 1.$
We have:
  $$\Vert \nabla^{t}_{\frac{\partial}{\partial s}}K^{t}_{i}(s)  \Vert_{\tilde{g}(t)}^{2} = (\dot{G})^{2}.$$
By the index lemma (e.g. \cite{Lee}), we deduce:
 $$I^{t}(J_{i}^{t},J_{i}^{t}) \leq I^{t}(K_{i}^{t},K_{i}^{t}), $$
and
 $$I^{t}(K_{i}^{t},K_{i}^{t})=\int_{0}^{r} \langle D_{s}K_{i}^{t},D_{s}K_{i}^{t} \rangle_{\tilde{g}(t)}-R_{m,\tilde{g}(t)} (K_{i}^{t},\dot{\gamma},\dot{\gamma},K_{i}^{t}) dt,$$
where $R_{m,\tilde{g}(t)} $ denote the $(4,0)$ curvature tensor associated to the metric $\tilde{g}(t)$. 
Hence:

\begin{eqnarray*}
\sum_{i=2}^{n}I^{t}(K_{i}^{t},K_{i}^{t}) &=& \sum_{i=2}^{n}\int_{0}^{r} \langle D_{s}K_{i}^{t},D_{s}K_{i}^{t} \rangle_{\tilde{g}(t)}-R_{m,\tilde{g}(t)} (K_{i}^{t},\dot{\gamma},\dot{\gamma},K_{i}^{t}) ds \\
 &=& \sum_{i=2}^{n} \int_{0}^{r} \Vert \nabla^{t}_{\frac{\partial}{\partial s}}K_{i}(s)  \Vert_{\tilde{g}(t)}^{2}-R_{m,\tilde{g}(t)} (K_{i}^{t},\dot{\gamma},\dot{\gamma},K_{i}^{t}) ds \\
 &=& \int_{0}^{r} (n-1)(\dot{G})^{2}- (G)^{2} Ric_{\tilde{g}(t)}(\dot{\gamma},\dot{\gamma})ds\\
 &\leq& (n-1)\int_{0}^{r} ((\dot{G})^{2} - (G)^{2} (\frac{1-\epsilon}{n-1})) ds .\\
\end{eqnarray*}
For performing the computation, we impose to $G$ to satisfy the O.D.E:
$$\left\lbrace \begin{array}{c}
   G(0) = G(r) = 1 \\
   \pp{G} + (\frac{1-\epsilon}{n-1}) G =0 \\
        \end{array} \right .$$
We notice that:
$$(\dot{G})^{2} - (G)^{2} (\frac{1-\epsilon}{n-1})) = (G\dot{G})', $$
and the solution of this O.D.E is given by the function:
$$G(s) = cos(\sqrt{\frac{1-\epsilon}{n-1}} s)+ \frac{1-cos(\sqrt{\frac{1-\epsilon}{n-1}} r)}{sin(\sqrt{\frac{1-\epsilon}{n-1}}r)} sin(\sqrt{\frac{1-\epsilon}{n-1}}s). $$
This function does not explode for  $ r $ in $[0 , \frac{\pi}{2\sqrt{\frac{1-\epsilon}{n-1}}} ], $
and, 
$$ (\dot G)(r) -(\dot G)(0) = -2\sqrt{\frac{1-\epsilon}{n-1}}\, tan(\frac{\sqrt{\frac{1-\epsilon}{n-1}}r}{2}).$$
Hence
$$ \sum_{i=2}^{n} I^{t}(J_{i}^{t},J_{i}^{t}) \leq -2 (n-1) \sqrt{\frac{1-\epsilon}{n-1}}\, tan(\frac{\sqrt{\frac{1-\epsilon}{n-1}}r}{2}) \leq 0 .$$
We get:
$$d\rho_{t}(Z^{1}_{t},Z^{3}_{t}) \leq C dt + 2d\beta_{t}. $$

After conditioning by $ \mathcal{F}_{T^{N}_{1}}$ we get the following computation:
$$\begin{array}{l}
 \P ( C_{N}  < (T^{N}_{1} + \frac12 )\wedge T^{N}_{2} ) \\
= \P (\exists t \in [(T^{N}_{1} ,(T^{N}_{1} + \frac12 )\wedge T^{N}_{2}]\text{ s.t.  }\rho_{t}(Z^{1}_{t},Z^{3}_{t})=0 )  \\
  \geq  \P \Bigg(\exists t \in [0, \frac12 ]\text{ s.t.  } Ct + 2\beta_{t} + \frac{ \frac{\pi}{\sqrt{1 +\epsilon}} \wedge \frac{C_{n}(d,K,cst-\epsilon) }{2}} {4} = 0 \\
 \text{ and  } \sup_{0\leq s \leq t} (Cs + 2\beta_{s} + \frac{ \frac{\pi}{\sqrt{1 +\epsilon}} \wedge \frac{C_{n}(d,K,cst-\epsilon) }{2}} {4}) < \frac{ \frac{\pi}{\sqrt{1 +\epsilon}} \wedge \frac{C_{n}(d,K,cst-\epsilon) }{2}} {2} \Bigg)\\
   \geq \tilde{\alpha} >0 .\\
 \end{array}$$
\end{proof}
\begin{remark}
A better $\tilde{\alpha}$ could be found with a martingale of the type $e^{a \beta_{t}-\frac{a^{2}}{2}t} $.
 \end{remark}

 \begin{stheorem}\label{unicite-g}
Let $(M,g)$ be a compact, strictly convex hypersurface isometrically embedded in $ \R^{n+1}$, $n \ge 2$, and $(M,g(t)) $ the family of metrics constructed by the mean curvature flow (as in \ref{prop-pullbac}). There exists a unique $g(T_{c}-t)$-BM in law.
\end{stheorem}

\begin{proof}
 Let $ X^{1}_{t}$ and $X^{2}_{t} $  two $g(T_{c}-t)$-BM , by a deterministic change of time we get two $\tilde{g}(t)$-BM that we note  $Z^{1}_{t}$ and $Z^{2}_{t}$. Let  $ N \leq T << 0$, as above we build $Z^{3}_{N,t} $, we obtain $Z^{3}_{N,t}=Z^{2}_{t} $ in law.
Let $\tilde{k}= E( T-N)$, where $E(t)$ is the integer part of $t$.
We have by construction:
$$\begin{array}{c}
\P ( \exists t\in [N, T] , \text{ s.t. } Z^{3}_{N,t} = Z^{1}_{t} ) \geq \P ( \exists t\in [T^{N}_{0}, T^{N}_{2\tilde{k}}], \text{ s.t. } Z^{3}_{N,t}= Z^{1}_{t}).
\end{array}$$
Let $\mathcal{F}$ be the natural filtration generated by the two processes, by propositions \ref{alfa1}, \ref{alfa2} and strong Markov property we obtain:
 $$\begin{array}{c}
 \P ( \exists t \in [N,T^{N}_{2}]  \text{ s.t. } Z^{3}_{N,t} = Z^{1}_{t} ) \\
\geq \P (T^{N}_{1} < \frac12 + N ; C_{N} < (T^{N}_{1} + \frac12 ) \wedge T^{N}_{2})\\
 = \E [ \P ( C_{N} \leq (T^{N}_{1} + \frac12 ) \wedge T^{N}_{2} | \mathcal{F}_{T^{N}_{1}}) \mathbbm{1}_{T^{N}_{1} \leq \frac12 + N} ] \\
 \geq  \tilde{\alpha} \E [ \mathbbm{1}_{T^{N}_{1} \leq \frac12 + N} ]\\
 \geq \alpha \tilde{\alpha} > 0 .\\
\end{array}$$
By successive conditioning (by $\mathcal{F}_{T_{2\tilde{k}-2}}$, ... ) we get:
$$ \P (\nexists t \in [T^{N}_{0}, T^{N}_{2\tilde{k}}] \text{ s.t. } Z^{3}_{N,t}= Z^{1}_{t}) \leq (1-\alpha \tilde{\alpha})^{\tilde{k}}.$$
Let $f_{1}...f_{m} \in \mathcal{B}_{b} (M)$ (bounded Borel functions ) and $ t<t_{1}<...<t_{m} \leq 0$,
 \begin{align*}
&  |\E [f_{1}(Z^{1}_{t_{1}})...f_{m}(Z^{1}_{t_{m}})-f_{1}(Z^{2}_{t_{1}})...f_{m}(Z^{2}_{t_{m}}) ] |  \\
&\qquad = |\E [f_{1}(Z^{1}_{t_{1}})...f_{m}(Z^{1}_{t_{m}})-f_{1}(Z^{3}_{N,t_{1}})...f_{m}(Z^{3}_{N,t_{m}})] |\\
&\qquad \leq \E [|f_{1}(Z^{1}_{t_{1}})...f_{m}(Z^{1}_{t_{m}})-f_{1}(Z^{3}_{N,t_{1}})...f_{m}(Z^{3}_{N,t_{m}})|\1_{Z^{1}_{t} \neq Z^{3}_{N,t}} ]\\
&\qquad \leq 2\Vert f_{1} \Vert_{\infty}...\Vert f_{m} \Vert_{\infty} \P(Z^{1}_{t} \neq Z^{3}_{N,t}) \\
&\qquad = 2\Vert f_{1} \Vert_{\infty}...\Vert f_{m} \Vert_{\infty} \P(\nexists u \in [N,t], \text{ s.t. } Z^{1}_{u} = Z^{3}_{N,u}) \\
&\qquad \leq  2\Vert f \Vert_{\infty}...\Vert f_{m} \Vert_{\infty} (1-\alpha \tilde{\alpha})^{E(t-N)}
 \end{align*}
 We get the result by sending $N$ to $-\infty$.
\end{proof}

 As application, we give uniqueness of a solution of a differential equation without initial condition.

\begin{corollary}\label{unicite-EDP}
Let $(M,g)$ be a compact, strictly convex hypersurface isometrically embedded in $ \R^{n+1}$, $n \ge 2$, and $(M,g(t)) $ the family of metrics constructed by the mean curvature flow (as in \ref{prop-pullbac}).
Then the following equation has a unique solution in $]0,T_{c}]$, where $T_{c}$ is the explosion time of the mean curvature flow.

\begin{equation}
\left\{ \begin{array}{l}\label{ed}
\frac{\partial}{\partial t}h(t,y) + H^2(T_{c}-t,y)h(t,y) = \frac12 \D_{g(T_{c}-t)}h(t,y)\\
\int_{M}h(T_{c},y)d\mu_{0}= 1 
 \end{array}\right.
\end{equation}
\end{corollary}
\begin{proof}
Existence:
let $X^{T_{c}}_{]0,T_{c}]}$ be a $g(T_{c}-t)$-BM with law at time $t$, $h(t,y)d\mu_{T_{c}-t}$. Then the law satisfies the equation (\ref{ed}), it is a consequence of a Green formula (compare with the similar computation for the Ricci flow in \cite{metric} section 2).
 
Uniqueness:
let $\tilde{h}$ be a solution of (\ref{ed}), and $\nu_{k}$ be a non-increasing sequence in $]0..T_{c}]$ such that $\lim_{k \rightarrow \infty} \nu_{k} =0$. Take a $M$-valued random variable  $\tilde{X}^{\nu_{k}} \sim \tilde{h}_{\nu_{k}} d\mu_{T_{c}-\nu_{k}}$, define the process:

\begin{displaymath}
 \overline{X}^{\nu_{k}}_{t}= \left\{ \begin{array}{ll}
  \tilde{X}^{\nu_{k}} & \textrm{ for   $t \in ]0..\nu_{k}]$ }\\
  g(T_{c}-t)\text{-BM} (\tilde{X}^{\nu_{k}}) & \textrm{ for  $  t \in [\nu_{k}..T_{c}]$ }
   \end{array} \right.
 \end{displaymath}

 By the similar argument as in section 2, we deduce the tightness of the sequence $ \overline{X}^{\nu_{k}}$, let $\overline{X}$ be a limit of a extracted sequence (also noted by $\nu_{k}$). 
 It is easy to see (by uniqueness of a solution of S.D.E, and of P.D.E with starting function) that $\overline{X}^{\nu_{k'}}_{(.)} \overset{\mathcal{L}}{=} \overline{X}^{\nu_{k}}_{(.)}$ for times greater than $ \nu_{k}$ and $k' \ge k $ . Sending $k'$ to infinity, we obtain $\overline{X}_{(.)} \overset{\mathcal{L}}{=} \overline{X}^{\nu_{k}}_{(.)}$ for times greater than $ \nu_{k}$. Note also that for $t \ge \nu_{k}$
$$\overline{X}^{\nu_{k}}_{(.)} \overset{\mathcal{L}}{=} g(T_{c}-.)\text{-BM}(\overline{X}^{\nu_{k}}_{t}) \overset{\mathcal{L}}{=} g(T_{c}-.)\text{-BM}(\overline{X}_{t}).$$ Hence $\overline{X}$ is a $g(T_{c}-t)_{]0,T_{c}]} $ Brownian motion. For $t\ge \nu_{k}$ we have  $$\overline{X}_{t} \overset{\mathcal{L}}{=} \overline{X}^{\nu_{k}}_{t} \sim \tilde{h}_{t}d\mu_{T_{c}-t} .$$
By uniqueness in law of such process, we get the uniqueness of the solution, hence $h = \tilde{h}$.
 \end{proof}
\bibliographystyle{plain}  
\bibliography{/home/kolehe/these_kolehe/biblio}
\end{document}